\documentstyle{article}
\begin{document}
\baselineskip+6pt \small \begin{center}{\textit{ In the name of
Allah, the Beneficent, the Merciful}}\end{center}
 \large
\begin{center}{Ergodicity of power series-map on the simplex of group algebra of a finite group.
} \end{center}
\begin{center} {Ural Bekbaev \footnote[1]{e-mail: bekbaev@science.upm.edu.my},
Mohamat Aidil M. J.\footnote[2]{e-mail: aidil@science.upm.edu.my}}
\end{center}
\begin{center} {Department of Mathematics $\&$ Institute for Mathematical
Research,}\end{center}
\begin{center} {FS,UPM, 43400, Serdang, Selangor, Malaysia.}\end{center}

\begin{abstract}{A finite group
$G$, its group algebra $R[G]$ over the field of real numbers, any
power series \\ $p(t)= a_0+a_1t+ a_{2}t^{2}+ ...$ , where $ a_i
\geq 0$, and $a_0+a_1+ a_{2}+...= 1$, and simplex
$$ S= \{x=\sum_{g\in G}x_gg\in R[G]: \sum_{g\in G}x_g=1, x_g\geq 0 \mbox{ for
any}\quad g\in G \}$$ are considered. Ergodicity  of the map $p:
S\rightarrow S$, where $p(x)= a_0+a_1x+ a_{2}x^{2}+ ...$ for $x\in
S$, on $S$ is shown. The regularity of this map at a given point
$x\in S $ is investigated as well.}\end{abstract} \vspace{1cm}

{\Large \bf 1. Introduction and preliminary}

The convergence problems of the sequences $\{p^{[n]}(x)\}_{n\in N
}$ (regularity of the map $p$ at $x$) and $\{q^{(n)}(x)\}_{n\in N
}$(ergodicity of the map $p$ at $x$) , where $p: S\rightarrow S$
is a given map, $S$ is a simplex in a finite dimensional real
vector space, $x\in S$, $p^{[1]}(x)=p(x)$,
$p^{[n+1]}(x)=p(p^{[n]}(x))$,
$q^{(n)}(x)=\frac{1}{n}\sum_{i=1}^np^{[i]}(x)$ for $n\in N$, are
important problems in applications of Mathematics in other areas
of Science (see, for example,[1],[2]). We study these problems
when $S$ is the natural simplex in the group algebra of a finite
group over the real numbers and $p$ is a power series. The
polynomial case of $p(t)$ has been considered in [3] and [4].

First let us prove some common results to use in the future. Let
$\{p^{[n]}(t)=\sum_{k=0}^{\infty}a_k^{[n]}t^k\}_{n\in N }$ be any
sequence of power series for which $\sum
_{k=0}^{\infty}{|a_k^{[n]}|}<\infty$ and
$\lim_{n\rightarrow\infty}a^{[n]}=0$, \\
where $a^{[n]}=\sup\{\vert{a_k^{[n]}}\vert : k\in N\}$. It should
be noted that the sequence of series $\{q^{(n)}(t)\}_{n\in N
}=\{\sum_{k=1}^{\infty}b_k^{(n)}t^k\}_{n\in N }$ also holds the
same property i.e. $\lim_{n\rightarrow\infty}b^{(n)}=0$, where
$b^{(n)}=\sup\{\vert{b_k^{(n)}}\vert : k\in N\}$ and
$q^{(n)}(t)=\frac{1}{n}\sum_{i=1}^np^{[i]}(t)$. Let $p_0^{[n]}(t)$
stand for $p^{[n]}(t)-a_0^{[n]}$.

Assume  that $A$ is a finite dimensional mono associative algebra
over $R$ i.e. every subalgebra of $A$ generated by one element is
an associative algebra. In future the vector space $A$ will be
considered only with respect to Euclidian topology.

 {\bf Lemma 1.} \textit{ If $x\in A$ is such an element that
$\lim_{n\rightarrow \infty}x^n=0 $ then
$$\lim_{n\rightarrow \infty}p_0^{[n]}(x)=0$$ as well.}

Proof of this lemma is not difficult due to the finite dimension
of $A$, assumption $\lim_{n\rightarrow\infty}a^{[n]}=0$ and
equality $\lim_{n\rightarrow \infty}x^n=0 $ or see [3].

{\bf Lemma 2.} \textit{ If $x\in A$ and $\lim_{n\rightarrow \infty
}x^n= x_0$ then $$\lim_{n\rightarrow \infty }p^{[n]}(x)$$ exists
if and only if $\lim_{n\rightarrow \infty
}[a_0^{[n]}e+(p^{[n]}(1)-a_0^{[n]})x_0]$ exists, where $e$ stands
for the identity element of $A$. Moreover when it is the case then
$$\lim_{n\rightarrow \infty }p^{[n]}(x)=\lim_{n\rightarrow \infty
}[a_0^{[n]}e+(p^{[n]}(1)-a_0^{[n]})x_0]$$}

{\bf Proof.} Due to $x_0 =\lim_{m\rightarrow \infty
}x^{m+1}=x\lim_{m\rightarrow \infty }x^{m}=(\lim_{m\rightarrow
\infty }x^{m})x$ and $x_0 =\lim_{m\rightarrow \infty
}x^{2m}=\lim_{m\rightarrow \infty }x^{m}\lim_{m\rightarrow \infty
}x^{m}$ one has $x_0 =xx_0 =x_0x =x_0^2$, in particular
$x_0(x-x_0)=(x-x_0)x_0= 0$.

Therefore $$x^n=(x_0+(x-x_0))^n=x_0^n+(x-x_0)^n= x_0+(x-x_0)^n$$
and $\lim_{n\rightarrow \infty }(x-x_0)^n= 0$.  Moreover
$$p^{[n]}(x)=a_0^{[n]}e+\sum_{i=1} ^{\infty }a_i^{[n]}(x_0+(x-x_0))^i=a_0^{[n]}e+\sum_{i=1} ^{\infty
}a_i^{[n]}(x_0+(x-x_0)^i)=$$
$$a_0^{[n]}e+p_0^{[n]}(1)x_0+p_0^{[n]}(x-x_0)=a_0^{[n]}e+p_0^{[n]}(1)x_0 +p_0^{[n]}(x-x_0)$$
and due to Lemma 1 one has $\lim_{n\rightarrow \infty
}p_0^{[n]}(x-x_0)= 0$. Therefore $\lim_{n\rightarrow \infty
}p^{[n]}(x)$ exists if and only if $\lim_{n\rightarrow \infty
}[a_0^{[n]}e+(p^{[n]}(1)-a_0^{[n]})x_0]$ exists.

Now we will consider any finite group $G$ and its group algebra
$A=R[G]$ over the field of real numbers [5]. For any $x\in R[G]$
the number $x_g$ stands for the coefficient at $g$ in the linear
expansion of $x$ with respect to the basis $\{g: g\in G\}.$

Let $L: R[G]\rightarrow R$ be map defined by $L(x)= \sum_{g\in
G}x_g$ for any $x=\sum_{g\in G}x_gg\in R[G]$ and $$S= \{x\in R[G]:
L(x)=1, x_g\geq 0 \mbox{ for any}\quad g\in G \}$$

It is clear that $S$ is a compact set and moreover it is closed
with respect to the multiplication as far as $L(xy)=L(x)L(y)$ for
any $x,y\in R[G]$.

Let $\mbox{Supp}(x)= \{g: x_g\neq 0\}$ and $c$ stand for
$\frac{1}{\vert G \vert}\sum_{g\in G}g\in \stackrel{\circ}{S}$,
where $$\stackrel{\circ}{S}=\{x\in S: x_g >0 \mbox{ for any} \quad
g\in G \}$$ It is evident that $cS=Sc=\{c\}.$

{\bf Proposition 1.} \textit{ If $y, z\in R[G]$ any two elements
with nonnegative components then \\1.
$\mbox{Supp}(y+z)=\mbox{Supp}(y)\cup \mbox{Supp}(z)$\\2.
$\mbox{Supp}(yz)=\mbox{Supp}(y)\mbox{Supp}(z),$ where
$\mbox{Supp}(y)\mbox{Supp}(z)$ stands for $\{gh: y\in
\mbox{Supp}(y), h\in \mbox{Supp}(z)\}$. In particular,
$\vert\mbox{Supp}(yz)\vert \geq
\mbox{max}\{\vert\mbox{Supp}(y)\vert, \vert\mbox{Supp}(z)\vert
\}$, and, if $y_e\neq 0 $,$z_e\neq 0 $ then $\mbox{Supp}(y)\cup
\mbox{Supp}(z)\subset \mbox{Supp}(yz)$}

The proof is evident.

In future $x$ will stand for a fixed element of $S$ and the
following notations will be used:
 $$n_x=\min\{ k\in N: (x^k)_e\neq 0\}=\min \{ k\in N: e\in \mbox{Supp}(x^k)\},$$
 where $e$ stands for the neutral element of $G$, $$G_x=\langle
\mbox{Supp}(x^{n_x})\rangle $$ is the subgroup of $G$ generated by
$\mbox{Supp}(x^{n_x})$, $c_x=\frac{1}{\vert G_x \vert}\sum_{g\in
G_x}g$ and $$ m_x=\min\{ k\in N: x^k\in R[G_x]\}=\min\{ k\in N:
\mbox{Supp}(x^k)\subset G_x\}$$ It is clear that $m_x$ divides
$n_x$ whenever $x\in S$.

{\bf Proposition 2.} \textit{ 1. The system of elements
$\{e,x,x^2,...,x^{n_x-1}\}$ is linear independent. \\ 2. If $r_1$,
$r_2$ are any nonnegative integers then $c_xx^{r_1}=c_xx^{r_2}$ if
and only if $r_1=r_2\mbox{mod} m_x$.}

{\bf Proof.} If the system $\{e,x,x^2,...,x^{n_x-1}\}$ is linear
dependent then $$x^m=
a_{m-1}x^{m-1}+a_{m-2}x^{m-2}+...+a_{1}x+a_0e$$ for some $m< n_x$
and real numbers $a_{m-1}, a_{m-2}, ... ,a_{0}$. Multiplication
the above equality by $x^{n_x-m}$ shows that $x^{n_x}$ can be
represented as a linear combination of $\{x,x^2,...,x^{n_x-1}\}$
which is a contradiction as far as $(x^{n_x})_e \neq 0$ and
$x_e=(x^2)_e=...=(x^{n_x-1})_e=0$. It is the proof of the first
part of Proposition 2.

If $n=km_x+r$, where $0\leq r < m_x$, then
$c_xx^n=c_x(x^{m_x})^kx^r=c_xx^r$ because of $x^{m_x}\in R[G_x]$,
$x\in S$. If $0\leq r_1 < r_2< m_x$ and $c_xx^{r_1}=c_xx^{r_2}$
then $$c_xx^{r_1+m_x-r_2}=c_xx^{r_2+m_x-r_2}=c_xx^{m_x}\in
R[G_x]$$ which evidently implies that $r_1+m_x-r_2< m_x$ and
$x^{r_1+m_x-r_2}\in R[G_x]$.This contradiction completes the proof
of the second part of Proposition 2.

{\bf Proposition 3.} \textit{ If $y\in S$ and $y_e\neq 0 $ then
the following conditions are equivalent.
$$ 1.\quad G_y=G $$
$$ 2.\quad y^n\in \stackrel{\circ}{S}\mbox{ for some} \quad n\in N $$}

{\bf Proof.} Due to Proposition 1 one has $\mbox{Supp}(y^n)\subset
\mbox{Supp}(y^m)$ whenever $n\leq m$. If $G_y=G $ then for every
$g\in G$ one can find such $n\in N$ and $g_1,g_2,...,g_n \in
\mbox{Supp}(y)$ for which $g=g_1g_2...g_n$. Therefore $(y^n)_g\neq
0$. Due to finiteness of $G$ it implies that $y^m\in
\stackrel{\circ}{S}$ for all big enough $m\in N$. The implication
of second condition the first is evident.

 Let $S(x)$ stand for the set of all limits of all converging subsequences of
$\{x^m\}_{m\in N}$.

{\bf Theorem 1.} \textit{ For any $x\in S$ the equalities
$$xc_x=c_xx,$$
$$S(x)= \{c_xx^r: 0\leq r< m_x\}$$
are valid. In particular, $\vert S(x)\vert=1$ i.e. $S(x)=\{c_x\}$
if and only if $x\in R[G_x]$.}

{\bf Proof.}  Consider the following expansion of the algebra $$
R[G]= Rc\oplus I$$, where $I= \{x\in R[G]: L(x)=0\}$. With respect
to this expansion one has the following   representations  $ S= c+
S_0, \qquad \stackrel{\circ}{S}= c+ \stackrel{\circ}{S}_0, $ where
$$S_0= \{x=\sum_{g\in G}x_gg\in I: -\frac{1}{\vert G \vert}\leq
x_g\leq 1- \frac{1}{\vert G \vert} \quad \mbox{ for any}\quad
g\in G \}\qquad{ and}$$ $$\stackrel{\circ}{S}_0= \{x=\sum_{g\in
G}x_gg\in I: -\frac{1}{\vert G \vert}< x_g < 1- \frac{1}{\vert G
\vert} \quad \mbox{ for any}\quad  g\in G \}.$$

First let us consider the case $x\in \stackrel{\circ}{S}$. In this
case, $G_x= G$ and $x=c+x_0$, where $x_0 \in
\stackrel{\circ}{S}_0$. It is clear that $$x^m= c+x_0^m$$ for any
$m\in N$. Due to compactness of $S$ there is a subsequence
$\{m_k\}_{k\in N}$ of $\{m\}_{m\in N}$ such that
$\{x^{m_k}\}_{k\in N}=\{c+x_0^{m_k}\}_{k\in N}$ converges. But for
the given $x_0$ one can find number $\lambda >1$ such that
$\lambda x_0 \in \stackrel{\circ}{S}_0$ so that
$\overline{x}=c+\lambda x_0 \in \stackrel{\circ}{S}.$ Once again
due to compactness of $S$ one can find subsequence
$\{m_{k_l}\}_{l\in N}$ of $\{m_{k}\}_{k\in N}$ such that the
sequence $\{\overline{x}^{m_{k_l}}\}_{l\in N}=
\{c+\lambda^{m_{k_l}}x_0^{m_{k_l}}\}_{l\in N}$ also converges. But
convergence of sequences $\{x^{m_{k_l}}\}_{l\in N}=
\{c+x_0^{m_{k_l}}\}_{l\in N}$, $\{\overline{x}^{m_{k_l}}\}_{l\in
N}= \{c+\lambda^{m_{k_l}}x_0^{m_{k_l}}\}_{l\in N}$ and the
inequality $\lambda >1$ imply that $$\lim_{k\rightarrow \infty
}x_0^{m_{k}}=0,\qquad   \lim_{k\rightarrow \infty }x^{m_{k}}=c.$$

Now assume that the sequence $\{x^m\}_{m\in N}$ does not converge
to $c$. In this case one can find a subsequence $\{x^{m_s}\}_{s\in
N}$ of it which itself has no subsequence converging to $c$. But
due to compactness of $S$ one can find a converging subsequence
$\{x^{m_{s_k}}\}_{k\in N}$ of $\{x^{m_s}\}_{s\in N}$. As it has
been shown already in such case $$\lim_{k\rightarrow \infty
}x^{m_{s_k}}=c.$$ This contradiction shows that
$\lim_{m\rightarrow \infty }x^{m}=c$ and completes the proof of
the theorem in the case $x\in \stackrel{\circ}{S}$.

Now let us consider case $x\in S\setminus \stackrel{\circ}{S}_0$
such that $x^n \in \stackrel{\circ}{S}_0$ for some fixed $n\in N$.
Due to the proof of the previous case $\lim_{m\rightarrow \infty
}x^{nm}=c.$ Let $\{m_k\}_{k\in N}$ be any subsequence of
$\{m\}_{m\in N}$, for which $\{x^{m_k}\}_{k\in N}$ converges, and
$m_k= nq_k+ r_k$, where $q_k$ is an integer and $0\leq r_k < n$
for any $k\in N$. In this case $\lim_{k\rightarrow \infty
}x^{m_k}= \lim_{k\rightarrow \infty }x^{nq_k+
r_k}=\lim_{k\rightarrow \infty }[(x^n)^{q_k}x^{r_k}]=
\lim_{k\rightarrow \infty }[x^{r_k}(x^n)^{q_k}]$. Without loss of
generality one can assume that $r_k = r$ for all $k\in N$.
Therefore $\lim_{k\rightarrow \infty }x^{m_k}= (\lim_{k\rightarrow
\infty }(x^n)^{q_k})x^r= cx^r =x^rc= c$ as far as $x=c+x_0$ and
$cx_0=x_0c=0$.

Now let $x$ be any element of $S$. One has $x^{m_x}\in R[G_x]$ and
due to Proposition 3 there is $n\in N$ for which
$$x^{nn_x}\in \stackrel{\circ}{S_{G_x}}=\{y=\sum_{g\in G_x}y_gg: L(y)=1
\quad\mbox{and}\quad y_g>0 \mbox{ for any} \quad g\in G_x \}.$$ It
implies that $x^{mm_x}\in \stackrel{\circ}{S_{G_x}}$ for some
$m\in N$ a far as $m_x$ divides $n_x$. Therefore due to the first
part of the proof
$$\lim_{k\rightarrow \infty }(x^{m_x})^k= c_x.$$

Assume that a sequence $\{x^{m_k}\}_{k\rightarrow \infty}$
converges and $m_k= m_xq_k+ r_k$, where $0\leq r_k< m_x$ for any
$k\in N$. One can assume that $r_k=r$ for all $k\in N$. Therefore
$$ \lim_{k\rightarrow \infty}x^{m_k}= \lim_{k\rightarrow
\infty}[(x^{m_x})^{q_k}x^r]=\lim_{k\rightarrow
\infty}[(x^{m_x})^{q_k}]x^r= x^r\lim_{k\rightarrow
\infty}[(x^{m_x})^{q_k}]=c_xx^r=x^rc_x$$ So  $$S(x)= \{c_xx^r:
0\leq r< m_x\}=\{x^rc_x: 0\leq r< m_x\}$$ and due to Proposition 2
one has $\vert S(x)\vert =1$ if and only if $m_x=1$ i.e. $x\in
R[G_x]$. This is the end of proof of Theorem 1.

{\bf Proposition 4.} \textit{ 1. For any $x\in S$ the inequality
$n_{xc_x}\leq m_x$ is true. \\ 2. If $n_{xc_x}=n_x$ then
$n_x=m_{x}=m_{xc_x}$}

{\bf Proof.} To show the inequality $n_{xc_x}\leq m_{x}$ assume
that $\mbox{Supp}(x^{m})\subset G_x$ for some $m\in N$.  In this
case $e\in \mbox{Supp}((xc_x)^{m})$ as far as $xc_x=c_xx$ due to
Theorem 1 and $e\in \mbox{Supp}(x)^{m}G_x=
\mbox{Supp}((xc_x)^{m})$ which implies that $ n_{xc_x}\leq m$ i.e.
$n_{xc_x}\leq m_x$.

Now let us show the equalities $n_x=m_{x}=m_{xc_x}$.  According to
the first part one has $n_{xc_x}\leq m_x \leq n_x$ and therefore
$n_{xc_x}= m_x= n_x$. Due to $n_{xc_x}= n_x$ one has
$G_{xc_x}=\langle \mbox{Supp}((xc_x)^{n_{xc_x}})\rangle =\langle
\mbox{Supp}(x^{n_{x}})\mbox{Supp}(c_x)\rangle=\langle
\mbox{Supp}(x^{n_{x}})G_x\rangle = G_x$ as far as $G_x=\langle
\mbox{Supp}(x^{n_{x}})\rangle $. The equality
$$\mbox{Supp}((xc_x)^{n})=\mbox{Supp}(x^{n})G_x$$
for any $n\in N$ implies that inclusion
$\mbox{Supp}((xc_x)^{n})\subset G_{xc_x}=G_x$ holds if and only if
$\mbox{Supp}(x^{n})\subset G_x$ and therefore $m_{x}= m_{xc_x}.$

{\bf Lemma 3.} \textit{For any $x\in S$ either both limits
$\lim_{n\rightarrow \infty }p^{[n]}(x)$, $\lim_{n\rightarrow
\infty }p^{[n]}(xc_x)$ do not exist or both of them exist and
$$\lim_{n\rightarrow \infty }p^{[n]}(x)=\lim_{n\rightarrow \infty
}p^{[n]}(xc_x)$$}

{\bf Proof.} For the given element $x$ one can represent $
p^{[n]}(t)$ in the  following form
$$ p^{[n]}(t)= \sum_{k=0}^{\infty}a_k^{[n]}t^k=\sum_{r=0}^{m_x-1}a_r^{[n]}t^r+\sum_{r=0}^{m_x-1}t^r\sum_{k=1}^{\infty}a_{km_x+r}^{[n]}t^{km_x}$$
Therefore assuming that $x^{m_x}=c_x+x_1$ one has
$$ p^{[n]}(x)=\sum_{r=0}^{m_x-1}a_r^{[n]}x^r+\sum_{r=0}^{m_x-1}x^r\sum_{k=1}^{\infty}a_{km_x+r}^{[n]}(c_x+x_1)^{k}=$$
$$a_0^{[n]}e+\sum_{r=1}^{m_x-1}a_r^{[n]}x^r+c_x\sum_{r=0}^{m_x-1}x^r\sum_{k=1}^{\infty}a_{km_x+r}^{[n]}+\sum_{r=0}^{m_x-1}x^r\sum_{k=1}^{\infty}a_{km_x+r}^{[n]}x_1^{k}=$$
$$c_x\sum_{r=0}^{m_x-1}x^r\sum_{k=0}^{\infty}a_{km_x+r}^{[n]}+a_0^{[n]}(e-c_x)+
(e-c_c)\sum_{r=1}^{m_x-1}a_r^{[n]}x^r+\sum_{r=0}^{m_x-1}x^r\sum_{k=1}^{\infty}a_{km_x+r}^{[n]}x_1^{k}$$
and
 $$ p^{[n]}(xc_x)=\sum_{r=0}^{m_x-1}a_r^{[n]}(xc_x)^r+\sum_{r=0}^{m_x-1}(xc_x)^r\sum_{k=1}^{\infty}a_{km_x+r}^{[n]}((c_x+x_1)c_x)^{k}=$$
$$a_0^{[n]}e+c_x\sum_{r=1}^{m_x-1}a_r^{[n]}x^r+c_x\sum_{r=0}^{m_x-1}x^r\sum_{k=1}^{\infty}a_{km_x+r}^{[n]}$$
$$=c_x\sum_{r=0}^{m_x-1}x^r\sum_{k=0}^{\infty}a_{km_x+r}^{[n]}+a_0^{[n]}(e-c_x)$$
 Due to Lemma 1 one has
$$\lim_{n\rightarrow \infty }\sum_{r=1}^{m_x-1}a_r^{[n]}x^r=\lim_{n\rightarrow \infty }\sum_{r=0}^{m_x-1}x^r\sum_{k=1}^{\infty}a_{km_x+r}^{[n]}x_1^{k}=0$$
as far as $\lim_{k\rightarrow \infty }x_1^k=0$ and
$\lim_{n\rightarrow \infty }a^{[n]}=0$.

Therefore either all three limits $$\lim_{n\rightarrow \infty
}p^{[n]}(x),\quad \lim_{n\rightarrow \infty }p^{[n]}(xc_x),\quad
\lim_{n\rightarrow \infty
}(c_x\sum_{r=0}^{m_x-1}x^r\sum_{k=0}^{\infty}a_{km_x+r}^{[n]}+a_0^{[n]}(e-c_x))$$
exist and are equal or none of them exists.

So due to this Lemma convergence of $\{p^{[n]}(x)\}_{n\in N}$ can
be reduced to convergence of $\{p^{[n]}(xc_x)\}_{n\in N}$. One
advantage of it is that $n_{xc_x}\le m_x$. Therefore after finite
number of reductions one can reduce convergence of
$\{p^{[n]}(t)\}_{n\in N}$ at some point of $S$ to its convergence
at a point $x\in S$ for which $n_{xc_x}= n_x$.

{\bf Theorem 2. } \textit{If $x\in S$, $c_x\neq e$ and $n_{xc_x}=
n_x$ then the limit $\lim_{n\rightarrow \infty }p^{[n]}(x)$ exists
if and only if $\lim_{n\rightarrow \infty }a_{0}^{[n]}$,
$\lim_{n\rightarrow \infty }\sum_{k=0}^{\infty}a_{km_x+r}^{[n]}$
exist for all $0\le r< m_x$. If $x\in S$, $c_x= e$ and $n_{xc_x}=
n_x$ then the limit $\lim_{n\rightarrow \infty }p^{[n]}(x)$ exists
if and only if  $\lim_{n\rightarrow \infty
}\sum_{k=0}^{\infty}a_{km_x+r}^{[n]}$ exists for all $0\le r<
m_x$. When it is the case then}
$$\lim_{n\rightarrow \infty }p^{[n]}(x)=
c_x\sum_{r=0}^{m_x-1}x^r\lim_{n\rightarrow \infty }\sum_{k=0}^{\infty}a_{km_x+r}^{[n]}
+(e-c_x)\lim_{n\rightarrow \infty }a_{0}^{[n]}$$

{\bf Proof.} Due to Proposition 4 one has $m_{x}= n_x= n_{xc_x}=
m_{xc_x}$. If  $c_x\neq e$ then the system
$$e-c_x,c_x,xc_x, x^2c_x,...,x^{n_x-1}c_x$$ is linear independent.
Indeed, if
$$x^mc_x=
a_{m-1}x^{m-1}c_x+a_{m-2}x^{m-2}c_x+...+a_{1}xc_x+a_0c_x+b(e-c_x)$$,
for some $m< n_x$ and real numbers $a_{m-1}, a_{m-2}, ...
,a_{0},b$, then multiplication it by $(e-c_x)$ implies that $b=0$
so $$x^mc_x=
a_{m-1}x^{m-1}c_x+a_{m-2}x^{m-2}c_x+...+a_{1}xc_x+a_0c_x$$ Now
multiplication it by $x^{n_x-m}$ shows that $x^{n_x}c_x$ can be
represented as a linear combination of
$\{xc_x,x^2c_x,...,x^{n_x-1}c_x\}$ which is a contradiction as far
as $(xc_x)^{n_x}_e=(x^{n_x}c_x)_e \neq 0$ and
$$(xc_x)_e=(x^2c_x)_e=...=(x^{n_x-1}c_x)_e= 0$$ In the case of
$c_x=e$ still one has linear independence of
$$c_x,xc_x, x^2c_x,...,x^{n_x-1}c_x$$

Due to the proof of Lemma 3 existence of $\lim_{n\rightarrow
\infty }p^{[n]}(x)$ is equivalent to the existence of
$\lim_{n\rightarrow \infty
}(c_x\sum_{r=0}^{m_x-1}x^r\sum_{k=0}^{\infty}a_{km_x+r}^{[n]}+(e-c_x)a_{0}^{[n]})$
and due to the above linear independence the existence of the last
limit is equivalent to the existence of $\lim_{n\rightarrow \infty
}a_{0}^{[n]}$, $\lim_{n\rightarrow \infty
}\sum_{k=1}^{\infty}a_{km_x+r}^{[n]}$ for all $0\le r< m_x$. It is
the proof of Theorem 2.

{\bf Lemma 4.} \textit{ If $a_1, a_2, a_3,...$ is a sequence of
nonnegative real numbers such that $a_1+ a_2+ a_3+...= 1$ then for
any $k\geq 2$ and $2 \le i\le k$ the following inequality
$$\sum_{j_1+j_2+...+j_i=k}a_{j_1}a_{j_2}...a_{j_i} \le (1-a_k)a$$
is valid, where $a=\sup\{a_1,a_2,a_3,...\}$ and the expression
$\sum_{j_1+j_2+...+j_i=k}$ stands for the summation taken over all
$(j_1,j_2,...,j_i)$ with natural entries for which
$j_1+j_2+...+j_i=k$ . In particular
$$\sum_{j_1+j_2+...+j_i=k}a_{j_1}a_{j_2}...a_{j_i} \le a$$ for any
$k\geq 1$ and $1\le i\le k$.}

{\bf Proof.} $\sum_{j_1+j_2+...+j_i=k}a_{j_1}a_{j_2}...a_{j_i}=$
$$\sum_{j_1=1}^{k-(i-1)}a_{j_1}\sum_{j_2=1}^{k-(i-2)-j_1}a_{j_2}\sum_{j_3=1}^{k-(i-3)-(j_1+j_2)}a_{j_3}...\sum_{j_{i-1}=1}^{k-1-(j_1+j_2+...+j_{i-2})}a_{j_{i-1}}a_{k-(j_1+j_2+...+j_{i-1})}$$
$$\leq a\sum_{j_1=1}^{k-(i-1)}a_{j_1}\sum_{j_2=1}^{k-(i-2)-j_1}a_{j_2}\sum_{j_3=1}^{k-(i-3)-(j_1+j_2)}a_{j_3}...\sum_{j_{i-1}=1}^{k-1-(j_1+j_2+...+j_{i-2}}a_{j_{i-1}}$$ Due to $a_1+ a_2+ a_3+... = 1$ one has $$\sum_{j_1=1}^{k-(i-1)}a_{j_1}\sum_{j_2=1}^{k-(i-2)-j_1}a_{j_2}\sum_{j_3=1}^{k-(i-3)-(j_1+j_2)}a_{j_3}...\sum_{j_{i-1}=1}^{k-1-(j_1+j_2+...+j_{i-2})}a_{j_{i-1}}$$
$$\leq \sum_{j_1=1}^{k-(i-1)}a_{j_1}= 1-\sum_{j=k-i+2}^\infty a_j\leq 1-a_k$$
and therefore $$\sum_{j_1+j_2+...+j_i=k}a_{j_1}a_{j_2}...a_{j_i}
\leq (1-a_k)a$$ As to inequality
$$\sum_{j_1+j_2+...+j_i=k}a_{j_1}a_{j_2}...a_{j_i} \leq a$$ for
any $k\geq 1$ and $1\leq i\leq k$ it is now evident. \vspace{1cm}

{\Large \bf 2. Ergodicity of power series-map}

Let now $p(t)=\sum_{k=0}^{\infty}a_kt^k= a_0+a_1t+ a_{2}t^{2}+
...$ be any power series, where $0\le a_i<1$, and $p(1) =a_0+a_1+
a_{2}+...= 1$. Consider its iterations
 $p^{[n+1]}(t)= p(p^{[n]}(t)) $, where $n\in N$, $p^{([1]}(t)=p(t) $.
 It is clear that coefficients of power series $$ p^{[n]}(t)=
\sum_{k=0}^{\infty}a_k^{[n]}t^k$$
 are nonnegative and $p^{[n]}(1) =\sum_{k=0}^{\infty}a_k^{[n]}= 1$.
 The following result is about the behavior of $a_0^{[n]}$ and $a^{[n]}=
Sup\{a_k^{[n]}: k\in N\}$.

{\bf Theorem 3.} \textit{The sequence $\{a_0^{[n]}\}_{n\in N}$ is
a monotone increasing sequence, $\{a^{[n]}\}$ is monotone
decreasing sequence and $\lim_{n\rightarrow \infty}a^{[n]}= 0 $}

{\bf Proof.}  Let $p^{[n]}(t)=a_0^{[n]}+p_0^{[n]}(t)$, where
$a_0^{[n]}$ is the constant term of $p^{[n]}(t)$. Due to
$$p^{[n+1]}(t)=p(p^{[n]}(t))=\sum_{i=0}^{\infty}a_i (a_0^{[n]}+p_0^{[n]}(t))^i$$
and the binomial formula, one has
$$p^{[n+1]}(t)= \sum_{i=0}^\infty a_i \sum_{j=0}^i
a_0^{[n]i-j}p_0^{[n]}(t)^j (^i _j)=\sum_{i=0}^\infty a_i
[a_0^{[n]i}+\sum_{j=1}^i a_0^{[n]i-j}p_0^{[n]}(t)^j (^i _j)]$$
$$=\sum_{i=0}^\infty a_i a_0^{[n]i}+\sum_{i=1}^\infty a_i
\sum_{j=1}^i a_0^{[n]i-j}p_0^{[n]}(t)^j (^i _j)=\sum_{i=0}^\infty
a_i a_0^{[n]i}+\sum_{j=1}^\infty (\sum_{i=j}^\infty a_i
a_0^{[n]i-j}\frac{i!}{(i-j)!})\frac {p_0^{[n]}(t)^j}{j!}$$
$$=p(a_0^{[n]})+\sum_{j=1}^\infty\frac{p^{(j)}(a_0^{[n]})}{j!}p_0^{[n]}(t)^j.$$
So $$a_0^{[n+1]}=p(a_0^{[n]})$$ and
$$a_k^{[n+1]}=\sum_{i=1}^\infty
\frac{p^{(i)}(a_0^{[n]})}{i!}\sum_{j_1+j_2+...+j_i=k}a_{j_1}^{[n]}a_{j_2}^{[n]}
...a_{j_i}^{[n]},\eqno{(1)}$$ where $p^{(i)}$ stands for $i$-th
derivative of $p(t)$ and the last sum is taken over all $(j_i,
j_2,...j_i)$ with positive integer entries, for which
$j_i+j_2+...+j_i=k$. The function $p(t)$ is a monotone increasing
function on [0,1) therefore the sequence $\{a_0^{[n]}\}_{n\in N}$
is a monotone increasing one due to $a_0^{[n+1]}=p(a_0^{[n]})$ and
one has equality $p(a)=a$, where
$a=\lim_{n\rightarrow\infty}a_0^{[n]}$. If $a=1$ then
$\lim_{n\rightarrow \infty}a^{[n]}= 0 $ since $a_k^{[n]}\le
1-a_0^{[n]}$ whenever $k\in N$, so in this case, the theorem is
true.

For $f(t)=p(t)-t$ one has $f(0)= a_0 \ge 0$, $f(1)=0$, $f(a)=0$,
$f'(0)= a_1-1 < 0$ and $f'(t)$ is a monotone increasing function
on [0,1). If $f'(t)\le 0$ for all $t<1$ then the function $y=f(t)$
is a positive, strictly decreasing function on $[0,1)$ and due to
t $f(a)=f(1)=0$ one has $\lim_{n\rightarrow\infty}a_0^{[n]}=a=1$
and the theorem is valid.

Now assume that $\lim_{t\rightarrow 1}f'(t)> 0$. In this case
there is an unique $0< b <1$ for which $f'(b)=0$ and $0\le a <b$.
It should be noted that in this case $f'(a)<0$ i.e. $p'(a)<1$.

For any $k\geq 2$ due to $(1)$ one has
$$a_k^{[n+1]}=p'(a_0^{[n]})a_{k}^{[n]}+\sum_{i=2}^\infty\frac{p^{(i)}(a_0^{[n]})}{i!}
\sum_{j_1+j_2+...+j_i=k}a_{j_1}^{[n]}a_{j_2}^{[n]}...a_{j_i}^{[n]}=$$
$$p'(a_0^{[n]})a_{k}^{[n]}+\sum_{i=2}^\infty
\frac{p^{(i)}(a_0^{[n]})}{i!}(1-a_{0}^{[n]})^i
\sum_{j_1+j_2+...+j_i=k}\frac{a_{j_1}^{[n]}}{1-a_{0}^{[n]}}\frac{a_{j_2}^{[n]}}{1-a_{0}^{[n]}}...
\frac{a_{j_i}^{[n]}}{1-a_{0}^{[n]}}$$ and therefore due to Lemma
4, one has
$$a_k^{[n+1]}\le
p'(a_0^{[n]})a_{k}^{[n]}+\sum_{i=2}^\infty\frac{p^{(i)}(a_0^{[n]})}{i!}(1-a_{0}^{[n]})^i
(1-\frac{a_{k}^{[n]}}{1-a_{0}^{[n]}})\sup_{i\geq1}\frac{a_{i}^{[n]}}{1-a_{0}^{[n]}}=$$
$$
p'(a_0^{[n]})a_{k}^{[n]}+(1-\frac{a_{k}^{[n]}}{1-a_{0}^{[n]}})\frac{a^{[n]}}{1-a_{0}^{[n]}}
\sum_{i=2}^\infty \frac{p^{(i)}(a_0^{[n]})}{i!}(1-a_{0}^{[n]})^i$$
But $$\sum_{i=2}^\infty
\frac{p^{(i)}(a_0^{[n]})}{i!}(1-a_{0}^{[n]})^i=
1-p(a_0^{[n]})-(1-a_{0}^{[n]})p'(a_0^{[n]})$$ in particular
$1-p(a_0^{[n]})-(1-a_{0}^{[n]})p'(a_0^{[n]})\geq 0$, as far as
\\$1=p(1)=p(a_0^{[n]}+1-a_0^{[n]})=\sum_{i=0}^\infty a_i(a_0^{[n]}+(1-a_0^{[n]}))^i
=\sum_{i=0}^\infty
a_i\sum_{j=0}^i(^i_j)(a_0^{[n]})^{i-j}(1-a_0^{[n]})^j
\\= \sum_{i=0}^\infty (\sum_{i=j}^\infty (^i_j)a_i
(a_0^{[n]})^{i-j})(1-a_0^{[n]})^j= \sum_{j=0}^\infty \frac
{p^{(j)}(a_0^{[n]})}{j!}(1-a_0^{[n]})^j$ and therefore
$$a_k^{[n+1]}\le
p'(a_0^{[n]})a_{k}^{[n]}+(1-\frac{a_{k}^{[n]}}{1-a_{0}^{[n]}})\frac{a^{[n]}}{1-a_{0}^{[n]}}
(1-p(a_0^{[n]})-(1-a_{0}^{[n]})p'(a_0^{[n]}))\eqno {(2)}$$ Due to
$(1)$ one has $$a_1^{[n+1]}= p'(a_0^{[n]})a_{1}^{[n]}\le
p'(a_0^{[n]})a_{1}^{[n]}+(1-\frac{a_{1}^{[n]}}{1-a_{0}^{[n]}})\frac{a^{[n]}}{1-a_{0}^{[n]}}
(1-p(a_0^{[n]})-(1-a_{0}^{[n]})p'(a_0^{[n]}))$$which implies that
inequality (2) is valid at $k=1$ as well.

 First let us show that $\{a^{[n]}\}_{n\in N}$  is a monotone decreasing sequence.
Indeed due to inequality
 (2) for any $k\in N$ one has
$$a_k^{[n+1]}\le
p'(a_0^{[n]})a^{[n]}+\frac{a^{[n]}}{1-a_{0}^{[n]}}(1-p(a_0^{[n]})-(1-a_{0}^{[n]})p'(a_0^{[n]}))=a^{[n]}\frac{1-p(a_0^{[n]})}{1-a_{0}^{[n]}}\le
a^{[n]} $$ as far as
$\frac{1-p(a_0^{[n]})}{1-a_{0}^{[n]}}=\frac{1-a_0^{[n+1]}}{1-a_{0}^{[n]}}\le
1$ and therefore $a^{[n+1]}\le a^{[n]}$. In particular there
exists
$$\lim_{n\rightarrow \infty}a^{[n]}=\lambda \geq 0$$

Now due to inequality (2) one can have $a_k^{[n+1]}\le
a_{k}^{[n]}(p'(a_0^{[n]})-$
$$\frac{\lambda}{(1-a_{0}^{[n]})^2}(1-p(a_0^{[n]})
-(1-a_{0}^{[n]})p'(a_0^{[n]})))+\frac{a^{[n]}}{1-a_{0}^{[n]}}(1-p(a_0^{[n]})-(1-a_{0}^{[n]})p'(a_0^{[n]}))\eqno{(3)}
 $$

Consider function $g(t)=p'(t)-\lambda
\frac{(1-p(t)-(1-t)p'(t))}{(1-t)^2}$. It was shown that
$a_0^{[n]}$ tends to $a$ from the left side. Therefore for all $n$
big enough either \\\quad \mbox{(Case
1)}$$p'[a_{0}^{[n]}]-\frac{\lambda}
{(1-a_{0}^{[n]})^2}(1-p(a_0^{[n]})-(1-a_{0}^{[n]})p'(a_0^{[n]}))\geq
0 \quad  \mbox{or} $$
\\\quad  \mbox{(Case
2)}$$p'(a_{0}^{[n]})-\frac{\lambda}{(1-a_{0}^{[n]})^2}(1-p(a_0^{[n]})-(1-a_{0}^{[n]})
p'(a_0^{[n]}))< 0 $$

{\bf Case 1.} In this case for $n$ big enough inequality (3)
implies that $a^{[n+1]}\le$
$$ a^{[n]}(p'(a_0^{[n]})-\frac{\lambda}{(1-a_{0}^{[n]})^2}(1-p(a_0^{[n]})
-(1-a_{0}^{[n]})p'(a_0^{[n]})))+\frac{a^{[n]}}{1-a_{0}^{[n]}}(1-p(a_0^{[n]})-(1-a_{0}^{[n]})p'(a_0^{[n]}))$$
and therefore by taking limit one has $$\lambda \le
\lambda(p'(a)-\frac{\lambda}{1-a}(1-p'(a)))+ \frac{\lambda}
{1-a}(1-p'(a))(1-a)=\lambda (1-\frac{\lambda}{1-a}(1-p'(a)))$$
Assumption $\lambda\neq 0$ leads to a contradiction
 $0\le -\frac{\lambda}{1-a}(1-p'(a))$ as far as $1-p'(a)>0$.

{\bf Case 2.} In this case for $n$ big enough inequality (3)
implies that
$$a^{[n+1]}\le
\frac{a^{[n]}}{1-a_{0}^{[n]}}(1-p(a_0^{[n]})-(1-a_{0}^{[n]})p'(a_0^{[n]}))$$
and therefore by taking limit one has $$\lambda \le
\lambda(1-p'(a))$$ Assumption $\lambda\neq 0$ leads to a
contradiction  $0\le -p'(a).$ This is the end of the proof of
Theorem .

{\bf Corollary 1.} \textit{Let $p(t)$ be as in the Theorem 3 and
$q^{(n)}{(t)}=\frac 1{n}\sum_{i=1}^np^{[i]}(t)=\sum_{i=0}^\infty
q^{(n)}_{i}{t^i}$. Then $\{q_{0}^{(n)}\}_{n=1}^{\infty}$ is a
monotone increasing sequence and
$\lim_{n\rightarrow\infty}q^{(n)}=0$, where }$$q^{(n)}=
Sup\{q_{i}^{(n)}:i\in N\}$$

{\bf Corollary 2.} \textit{Let $p(t)=\sum_{k=0}^\infty a_k t^k$ be
any series such that $-1<a_i<1$ whenever $i\in N\cup \{0\}$and
$\sum_{i=0}^\infty {|a_i|}\leq1$ then $\lim_{n\rightarrow
\infty}a^{[n]}= 0 $, where $a^{[n]}= Sup\{|a_k^{[n]}|: k\in N\}$,
$p^{(1)}(t)=p(t) $, $p^{[n+1]}(t)= p(p^{[n]}(t)) $ and
$p^{[n]}(t)= \sum_{k=0}^{\infty}a_k^{[n]}t^k$}.

{\bf Corollary 3.} \textit{Let $p(t)$ be same as in Corollary 2.
Then if $\lim_{n\rightarrow \infty}p^{[n]}(0)=\lim_{n\rightarrow
\infty}a_0^{[n]}=a$ exists then $\lim_{n\rightarrow
\infty}p^{[n]}(t)=a$ for any $t\in (-1, 1)$}.

Theorem 2 motivates investigation of the following problem. Let
$p(t)=\sum_{i=0}^{\infty}a_{i}t^i$, where $0\leq a_i <1$ whenever
$i\in N\cup\{0\}$ and $p(1)=1$, be a power series, $m$ be a fixed
natural number. When does
$$\lim_{n\rightarrow \infty }\sum_{k=0}^{\infty}a_{km+r}^{[n]}$$
exist for all $0\le r< m$? Obviously it is equivalent to the
question: When does $\lim_{n\rightarrow \infty
}p^{[n]}(\overline{t})$ exist, where
$\overline{t}=t\mbox{mod}(t^m-1)$? Let
$\mbox{Supp}(p(t))=\mbox{Supp}(\sum_{i=0}^{\infty}a_{i}t^i)=
\{t^i: a_i \neq 0\ , i\in N \cup\{0\}\}$.

In future it is assumed that $p(t)= t^rp_0(t)$, where $r \geq 0$,
$\mbox{Supp} p_0(t)=\{ t^{q_i}: i<l\}$, where $2\leq l\leq \infty$
and $0=q_0< q_1< q_2<...$. It is not difficult to see that
$$p^{[n]}(t)=t^{r^n}p_0^{(n)}(t)\quad \mbox{and} \quad
p_0^{(n+1)}(t)=p_0^{(n)}(t)^rp_0(t^{r^n}p_0^{(n)}(t))$$
, which, due to Proposition 1, implies that \\
1. $\mbox{Supp}(p_0^{(n)}(t))\subset \langle
\{t^{q_i}:i<l\}\rangle$, where $\langle \{t^{q_i}: i<l\}\rangle$
stands for the semigroup generated by $\{t^{q_i}: i<l\}$.\\
2. $t^{q_i}\mbox{Supp}(p_0^{(n)}(t))\subset
\mbox{Supp}(p_0^{(n+1)}(t))$ for any $i<l$

The second property implies, as far as
$\mbox{Supp}(p_0^{(n)}(t))\subset \mbox{Supp}(p_0^{(n+1)}(t))$ ,
that  the sets $\mbox{Supp}(p_0^{(n)}(\overline{t}))$ are the same
sets for all big enough $ n\in N$. Let us denote it by
$G_{p_0(t)}$. Moreover $\overline{t}^{q_i}G_{p_0(t)}\subset
G_{p_0(t)}$ and $\{\overline{t}^{q_i}: i<l\}\subset G_{p_0(t)}$.
Therefore taking into account the first property one can conclude
that $G_{p_0(t)}$ is the subgroup of $\{\overline{t}^0,
\overline{t}^1,..., \overline{t}^{m-1}\}$ generated by
$\{\overline{t}^{q_i}: i<l\}$.
 Assume that it is as a cyclic group
generated by $\overline{t}^q$, where $0\leq q < m$.

It is clear that the sequence $\{r^k\mbox{mod} m\}_{k\in N}$ is a
repeating sequence i.e. there are $d\in N$ and different numbers
$m_0, m_1,...,m_{d-1}$ between $0$ and $m-1$ for which
$$r^{kd+i}= m_i\mbox{mod} m$$ whenever $i\in \{0,1,...,d-1\}$ and $k$ is big enough.

{\bf Theorem 4.} \textit{1. If
$\cap_{i=0}^{d-1}\overline{t}^{m_i}G_{p_0(t)}= \emptyset$ then
$\lim_{n\rightarrow \infty }p^{(n)}(\overline{t})$ does not exist.
\\2. If $\cap_{i=0}^{d-1}\overline{t}^{m_i}G_{p_0(t)}\neq
\emptyset$ then
$\overline{t}^{m_i}G_{p_0(t)}=\overline{t}^{m_j}G_{p_0(t)}$
whenever $i,j\in \{0,1,...,d-1\}$, the limit $\lim_{n\rightarrow
\infty }p^{[n]}(\overline{t})$ exists and $$\lim_{n\rightarrow
\infty
}p^{[n]}(\overline{t})=\overline{t}^{m_0}c_{p_0(t)}+(e-c_{p_0(t)})a$$,
where $c_{p_0(t)}=\frac{1}{\vert G_{p_0(t)} \vert}\sum_{g\in
G_{p_0(t)}}g$ and $a=\lim_{n\rightarrow \infty }a_0^{[n]}$}

{\bf Proof.} Assume that
$\cap_{i=0}^{d-1}\overline{t}^{m_i}G_{p_0(t)}= \emptyset$, $
\lim_{n\rightarrow \infty }p^{[n]}(\overline{t})$ exists  and the
coefficient at $\overline{t}^{m_0+h}$ in $ \lim_{k\rightarrow
\infty }p^{[kd]}(\overline{t})= \lim_{n\rightarrow \infty
}p^{[n]}(\overline{t})$ is not zero. Due to
$\cap_{i=0}^{d-1}\overline{t}^{m_i}G_{p_0(t)}= \emptyset $ one can
find $0\leq j<m$ for which $\overline{t}^{m_0+h}\notin
\overline{t}^{m_j}G_{p_0(t)}$. It means that if even
$\lim_{k\rightarrow \infty }p^{[kd+j]}(\overline{t})$ exists the
coefficient at $\overline{t}^{m_0+h}$ in $\lim_{n\rightarrow
\infty }p^{[n]}(\overline{t})$ is zero. This contradiction
indicates that $ \lim_{n\rightarrow \infty }p^{[n]}(\overline{t})$
can not exist.

Now assume that $\cap_{i=0}^{d-1}\overline{t}^{m_i}G_{p_0(t)}\neq
\emptyset.$ In this case one can find integer numbers $ s_0,
s_1,...,s_{d-1}$ such that
$\overline{t}^{m_0}\overline{t}^{qs_0}=\overline{t}^{m_1}\overline{t}^{qs_1}=...
=\overline{t}^{m_{d-1}}\overline{t}^{qs_{d-1}}$ and therefore
$\overline{t}^{m_i}G_{p_0(t)}=\overline{t}^{m_j}G_{p_0(t)}$
whenever $i,j\in \{0,1,...,d-1\}$. In particular
$\overline{t}^{r^n}c_{p_0(t)}=\overline{t}^{m_0}c_{p_0(t)}$ for
all $n\in N$ big enough.

We know that for $n_0\in N$ big enough $
p_0^{(n_0)}(\overline{t})\in  \stackrel{\circ}{S}_{G_{p_0(t)}}$
and due to Theorem 1 one has $\lim_{k\rightarrow \infty
}p_0^{(n_0)}(\overline{t})^k=c_{p_0(t)}. $ Let
$p_0^{(n_0)}(\overline{t})=c_{p_0(t)}+ x_1$. Due to
$\overline{t}^{r^{n_0}}p_0^{(n_0)}(\overline{t})=\overline{t}^{r^{n_0}}c_{p_0(t)}+\overline{t}^{r^{n_0}}x_1$
and $c_{p_0(t)}x_1=0$ one has
$$p^{[n]}(\overline{t}^{r^{n_0}}p_0^{(n_0)}(\overline{t}))=
p^{[n]}(\overline{t}^{r^{n_0}}c_{p_0(t)}+\overline{t}^{r^{n_0}}x_1)=
p^{[n]}(\overline{t}^{r^{n_0}})c_{p_0(t)}-a_0^{[n]}c_{p_0(t)}+
[p^{[n]}(\overline{t}^{r^{n_0}}x_1)-a_0^{[n]}e]+a_0^{[n]}e$$ But
$\lim_{k\rightarrow \infty }x_1^k=\lim_{k\rightarrow \infty
}(\overline{t}^{r^{n_0}}x_1)^k=0$ and therefore due to Lemma 2,
Theorem 3 one has $\lim_{n\rightarrow \infty }
p^{[n]}(\overline{t}^{r^{n_0}}x_1)-a_0^{[n]}e=0$. Therefore
$\lim_{n\rightarrow \infty
}p^{[n]}(\overline{t})=\lim_{n\rightarrow \infty
}p^{[n]}[\overline{t}^{r^{n_0}}p_0^{(n_0)}(\overline{t})]$ exists
if and only if $\lim_{n\rightarrow \infty
}(p^{[n]}(\overline{t}^{r^{n_0}})c_{p_0(t)}-a_0^{[n]}c_{p_0(t)}+
a_0^{[n]}e)$ exists. Now to finish the proof notice that for any
$n\in N$ one has $p_0^{(n)}(\overline{t}^{r^{n_0}})c_{p_0(t)} =
c_{p_0(t)}$ and
$p^{[n]}(\overline{t}^{r^{n_0}})=\overline{t}^{r^{n_0+n}}p_0^{(n)}(\overline{t}^{r^{n_0}})$
which implies that $\lim_{n\rightarrow \infty
}p^{[n]}(\overline{t})=\overline{t}^{m_0}c_{p_0(t)}+(e-c_{p_0(t)})a$,
where $a=\lim_{n\rightarrow \infty }a_0^{[n]}$.

{\bf Examples.} 1. Let $p(t)=\frac{1}{2}t^3(1+t^4)$, $m=12$. In
this case $r=3$, $p_0(t)=\frac{1}{2}(1+t^4)$,
$\overline{t}^{12}=1$, $G_{p_0(t)}=\{1, \overline{t}^4,
\overline{t}^8\}$, $d=2$ and $m_0=3$, $m_1=9$. The intersection
$\cap_{i=0}^{d-1}\overline{t}^{m_i}G_{p_0(t)}$ is empty as far as
$m_0G_{p_0(t)}=\{\overline{t}^3, \overline{t}^7,
\overline{t}^{11}\}$ and $m_1G_{p_0(t)}=\{\overline{t},
\overline{t}^5, \overline{t}^{9}\}$.

2. Now let $p(t)=\frac{1}{2}t^3(1+t^4)$ be same, $m=10$. In this
case $r=3$, $p_0(t)=\frac{1}{2}(1+t^4)$, $\overline{t}^{10}=1$,
$G_{p_0(t)}=\{1,\overline{t}^2, \overline{t}^4, \overline{t}^6,
\overline{t}^8,\}$, $d=4$ and $m_0=3$, $m_1=9$, $m_2=7$, $m_3=1$.
The intersection $\cap_{i=0}^{d-1}\overline{t}^{m_i}G_{p_0(t)}$ is
equal to $\{\overline{t}, \overline{t}^3, \overline{t}^5,
\overline{t}^7, \overline{t}^{9}\}$.

 {\bf Theorem 5.} \textit{The
above considered map $p:S\rightarrow S$ is ergodic on $S$.}

{\bf Proof.} Let $q^{(n)}(t)=\frac{1}{n}\sum_{i=1}^np^{[i]}(t)$
for $n\in N$. Due to Lemma 3,Theorem 2 and Theorem 3 to show
ergodicity of $p$ at $x\in S$ one can assume that $n_x= n_{xc_x}$
and prove existence of $$\lim_{n\rightarrow \infty
}\sum_{k=0}^{\infty}q_{km_x+r}^{(n)}$$ for all $0\leq r< m_x$. We
consider any fixed $m\in N$ and show that $\lim_{n\rightarrow
\infty }q^{(n)}(\overline{t})$ exists. It is easy to see that $
q^{(n)}[p^{[n_0]}(t)]=q^{(n+n_0)}(t)+\frac{n_0}{n}(q^{(n+n_0)}(t)-q^{(n_0)}(t))$
and \\ $\lim_{n\rightarrow \infty
}\frac{n_0}{n}(q^{(n+n_0)}(\overline{t})-q^{(n_0)}(\overline{t}))=0$.
So if  $\lim_{n\rightarrow \infty
}q^{(n)}(p^{[n_0]}(\overline{t}))$ exists then
$$\lim_{n\rightarrow \infty
}q^{(n)}[(p^{[n_0]}(\overline{t}))=\lim_{n\rightarrow \infty
}q^{(n)}(\overline{t})$$

For $n_0\in N$ big enough $ p_0^{(n_0)}(\overline{t})\in
\stackrel{\circ}{S}_{G_{p_0(t)}}$ and due to Theorem 1 one has
\\$\lim_{k\rightarrow \infty
}p_0^{(n_0)}(\overline{t})^k=c_{p_0(t)}. $ Let
$p_0^{(n_0)}(\overline{t})=c_{p_0(t)}+ x_1$. One has
$$q^{(n)}(p^{([n_0]}(\overline{t}))=q^{(n)}(\overline{t}^{r^{n_0}}p_0^{(n_0)}(\overline{t}))=
q^{(n)}(\overline{t}^{r^{n_0}}c_{p_0(t)}+\overline{t}^{r^{n_0}}x_1)=$$
$$q^{(n)}(\overline{t}^{r^{n_0}})c_{p_0(t)}-
q^{(n)}_{0}c_{p_0(t)}+[q^{(n)}(\overline{t}^{r^{n_0}}x_1)-q^{(n)}_{0}e]+
q^{(n)}_{0}e$$ and $\lim_{n\rightarrow \infty }
[q^{(n)}(\overline{t}^{r^{n_0}}x_1)-q^{(n)}_{0}]=0$ due to
$\lim_{k\rightarrow \infty }x_1^k=\lim_{k\rightarrow \infty
}(\overline{t}^{r^{n_0}}x_1)^k=0$. Moreover
$q^{(n)}(\overline{t}^{r^{n_0}})c_{p_0(t)}=
\frac{1}{n}\sum_{i=1}^n\overline{t}^{r^{n_0+i}}p_0^{(i)}[(\overline{t}^{r^{n_0}})c_{p_0(t)}$
and  $p_0^{(i)}(\overline{t}^{r^{n_0}})c_{p_0(t)}=c_{p_0(t)}$
whenever $i\in N$. Therefore
$$\lim_{n\rightarrow \infty }q^{(n)}(\overline{t})=
\lim_{n\rightarrow \infty }q^{(n)}[p^{(n_0)}(\overline{t})]=
\frac{1}{d}\sum_{i=0}^{d-1}\overline{t}^{m_i}c_{p_0(t)}+a(e-c_{p_0(t)})$$

To complete the picture let us consider the extremal case
$p(t)=t^r$ as well, where \\$r>1$. Let $S_r(x)$ stand for the set
of all limits all converging subsequences of $\{p^{(n)}(x)\}_{n\in
N}=\{x^{r^n}\}_{n\in N}$. The sequence $\{r^{n}\mbox{mod}
m_x\}_{n\in N}$ is a repeating sequence. Let
$(m_o,m_1,...,m_{d-1})$ be its repeating part, where
$m_0,m_1,...,m_{d-1}$ are different numbers between $0$ and
$m_x-1$. By the use of Theorem 1 the following result can be
proved easily.

{\bf Theorem 1'.} \textit{1. For any $x\in S$ the equalities
$$S_r(x)= \{c_xx^{m_i}: 0\leq i< d\}=\{x^{m_i}c_x: 0\leq i< d\}$$
are valid. In particular, $\vert S_r(x)\vert=1$ if and only if
$m_x$ divides $r^{k}(r-1)$ for some $k\in N$.\\ 2.
$\lim_{n\rightarrow \infty }q^{(n)}(x)=\lim_{n\rightarrow \infty
}\frac{1}{n}\sum_{i=0}^nx^{r^n}=
\frac{1}{d}\sum_{i=0}^{d-1}x^{m_i}c_x$ }\vspace{1cm}

{\Large\bf References}

(1) H. Kesten, Quadratic transformations: a model for population
growth. I,II. Adv. Appl. Probab. {\bf 2}(1970), 1-82; 179-228.

[2] Yu.I. Lyubich, Mathematical structures in population genetics,
Springer-Verlag, Berlin, 1992.

[3] U. Bekbaev,Regularity and ergodicity properties of polynomial
maps  induced by the multiplication of group algebra of a finite
group. Proceedings of The second International Conference on
Research and Education in Mathematics (ICREM 2),pp.134-144,May
25-27, 2005.

[4]  U. Bekbaev, Mohamat Aidil M. J.,Regularity and ergodicity
properties of polynomial maps induced by the multiplication of
group algebra of a finite group ( Nonzero constant case).
Proceedings of The second International Conference on Research and
Education in Mathematics (ICREM 2), pp.144-153,May 25-27, 2005.

[5] S.Lang, Algebra, Columbia University,Addison-Wesley,1972.

\end{document}